# Chapter 5

# Working times in atypical forms of employment: the special case of part-time work

Patrick LETREMY, Marie COTTRELL
*Samos-Matisse, CNRS UMR 8595, Université Paris 1, pley,cottrell@univ-Paris1. fr*

**Abstract**: In the present article, we attempt to devise a typology of forms of part-time employment by applying a widely used neuronal methodology called Kohonen maps. Starting out with data that we describe using category-specific variables, we show how it is possible to represent observations and the modalities of the variables that define them simultaneously, on a single map. This allows us to ascertain, and to try to describe, the main categories of part-time employment.

**Key words**: Kohonen maps, working times, classification.

## INTRODUCTION

France's economic recovery since 1997 has been accompanied by strong job creation and by a significant drop in unemployment. This does not mean however that there has been any real reduction in the number of people working under what has come to be known as "atypical forms of employment". Quite the contrary, the number of persons in temporary employment (with fixed term contracts hereafter FTC or doing temporary agency work) has never been as high. There has been an unprecedented rise in part-time work in France, something that coincides nowadays with the ever-increasing number of female entrants into the job market. In an Employment Survey carried out by the French National Statistics Office (INSEE), part-time jobs represented 16.8% of the country's employed working population, and temporary jobs (temporary work and FTC) 6.3%.

Furthermore, since 1994, there has been greater growth in temporary work than in FTC.

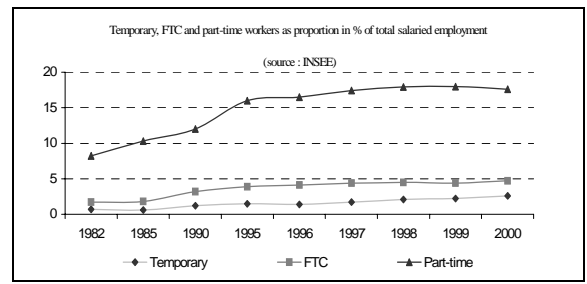

*Figure 1.* Changes in number of workers involved in atypical forms of employment

Such atypical forms of employment still constitute a relatively small minority of all jobs. It should be noted however that the circumstances surrounding female work have been considerably altered by the large and increasing proportion of women part-time workers. In addition, temporary work has a much greater effect on the labour market than the sheer weight of the numbers involved, a prime example being the preponderant role it plays in workers' transition back and forth between employment and unemployment.

The rise of these atypical forms of employment has unsurprising drawn attention to issues relating to the extent to which full-time employment will in the future be carried out by people working under an open-ended contract (*hereafter OEC*). The *IRES* and *MATISSE* research centres' contributions to the INSEE's 1998-99 *Timetables* survey focused on the working times which characterise these atypical forms of employment. For example, a study that benefited from a DARES research grant (Cottrell, Letrémy, Macaire et al. 2001), "*Working times with atypical forms of employment, Final Report*", IRES, Noisy-le-Grand, February 2001) tried to discover whether atypical forms of employment are subject to specific constraints in terms of the working times they entail. In other words, do they imply circumstances that should *a priori* be considered to be more difficult for those who are actually in this sort of work situation? Answering this question means making a comparison with a benchmark norm, the obvious one being the current situation for people in full-time employment and working under an open-ended contract. This means that we have not only tried to discover whether such atypical forms of employment are subject to specific time constraints, but also whether they are having to cope with working time constraints that are harder to deal with than is the case when the person involved benefits from an open-ended contract and a full-time employment status.

We should specify the terms which the present article uses:
– By atypical forms of employment we primarily mean i) the various modalities of temporary salaried work, whether full-time or part-time, and ii) the various modalities of part-time salaried work, regardless of the nature of the employment contract.
– By working times, not only do we mean issues relating to the number of hours worked, but also schedules, calendars and working times patterns, the variability and predictability thereof, how much choice the employee has in these different areas, etc. Note that the study only focuses on people's current principal activity.

Neuronal techniques such as Kohonen maps were used throughout the study to segment groups of employees according to available quantitative variables, before linking the category variable that is defined in this manner with informed qualitative variables. We would like use the present article to present an alternative to this technique, proposing a method that makes it possible to segment individuals by qualitative variables, even though this particular segmentation will later be crossed with available quantitative variables.

To present this new methodology, we took a particular interest in ***part-time employees*** working on either an open-ended or a fixed term contract. It is common knowledge that practically all part-time employment involves women (90% of OEC part-timers, 82.5% of FTC part-timers). This basically relates to women employees in areas such as retail, services and the social and non-profit sectors. However, we still wonder whether there are any differences between OEC and FTC part-timers – for example, whether the women who find themselves in either of these two situations have the same profile, whether they chose their part-time status or not, etc.?

We extracted data relating to part-time employees from from the INSEE's 1998-1999 *Timetables* survey. This covered 690 OEC and 137 FTC workers, after eliminating data that contained input errors or missing information. We then restricted the number of variable and kept 14 qualitative ones (type of contract, gender, age, the regularity of the timetabling, whether this involved sociable hours [night or weekend shifts], employee autonomy, the schedules predictability, etc), all of this for a total of 39 modalities. We also kept 5 quantitative variable relating to the number of hours worked per week. Data presentation takes the form of a complete disjunctive table containing 827 rows, 39 columns featuring 1s or 0s and 5 columns of real data (see the appendix for additional information on the survey).

*Table 1.* Qualitative variables

| Heading | Modalities | Name |
|---|---|---|
| Type of employment contract | Open-ended / fixed term contract | OEC, FTC |
| Gender | Man, Woman | MAN, FEM |
| Age | <25, [25, 40], [40,50], ≥50 | AGE1, AGE2, AGE3, AGE4 |
| Daily working schedules | Identical, Posted, Variable | HORIDE, HORP0S, HORVAR |
| Number of days worked per week | Identical, Variable | JWK1, JWK2 |
| Night shifts | Never, Sometimes, Usually | NITE1, NITE2, NITE3 |
| Saturday shifts | Never, Sometimes, Usually | SAT1, SAT2, SAT3 |
| Sunday shifts | Never, Sometimes, Usually | SUN1, SUN2, SUN3 |
| Wednesday shifts | Never, Sometimes, Usually | WED1, WED2, WED3 |
| Able to take time off | Yes, Yes under certain conditions, No | ABS1, ABS2, ABS3 |
| The schedule is determined by… | The firm itself, choice is given, he/her decides him/herself, other | DET1, DET2, DET3, DET4 |
| Part-time status forced | Yes, No | INVOL, VOL |
| Worker knows his/her schedule for next day(s) | Yes, No | LEND1, LEND2 |
| Possibility of carrying over working hours | Not applicable, Yes, No | RECUP0, RECUP1, RECUP2 |

*Table 2.* Quantitative variables

| Heading | Name |
|---|---|
| Minimum duration of actual workweek | DMIN |
| Maximum duration of actual workweek | DMAX |
| Theoretical duration of workweek | DTHEO |
| Number of overtime hours worked per week | HSUP |
| Number of hours of extended shift work/week | HPROL |

A simple cross-analysis of the variables reveals right away that men only represent 10% of all part-time employees working on an OEC basis and 18% of all part-time employees on an FTC. Moreover, even though forced (and therefore involuntary) part-time work accounts for nearly 50% of all employment contracts, it only represents 43% of the OEC, versus nearly 80% of FTC. Note that 83% of all contracts are OEC.

After a cursive study of these descriptive statistics (we will not be delving any further into them at present; see appendix for elements thereof), we are now going to carry out a segmentation of those individuals who are represented by the 14 qualitative variable defined above, as well as their 39 modalities. Towards this end, we will be defining a new method, one that is based on the Kohonen algorithm, but which enables an analysis of complete disjunctive tables.

## 1. THE KOHONEN ALGORITHM

This is the original classification algorithm that Teuvo Kohonen defined in the 1980s based on his studies of neuromimetic motivations (Kohonen 1984; 1995). In the present data analysis framework (Kaski 1997; Cottrell, Rousset 1997), the data space is a finite set that is identified by the rows of a data table. Each row of this table represents one of *N* individuals (or observations) that are being described by an identifier and by *p* quantitative variables. The algorithm then regroups the observations into separate classes, whilst respecting the topology of the data space.

This means that a priori we have defined a concept that accounts for a neighbourhood between classes. It also means that neighbouring observations in the data space of dimension *p* will belong (once they have been classified) to the same class or to neighbouring classes. The use of this algorithm is justified by the fact that it enables a regrouping of individuals into small classes whose neighbourhood is meaningful (unlike a hierarchical classification or a moving centres algorithm), and that they themselves can then be dynamically regrouped into super classes, preserving all the while the relationships of neighbourhood that have been detected. The visual representation of the classes is therefore easy to interpret, inasmuch as it occurs at a global level. Inversely, visual representations obtained through the use of classical projection methods are incomplete, as it becomes necessary to consult a number of successive projections in order to derive any reliable conclusions.

The structures of neighbourhood between the various classes can be chosen in a variety of ways, but in general we assume that the classes are laid out on a rectangular two-dimensional grid, this being a natural definition of neighbours in each class. We can also consider a one-dimensional topology, a so-called string, and possibly even a toroidal structure or a cylinder.

## 1.1     The algorithm for the quantitative data

The classification algorithm is an iterative one. It is launched through the association of each class with a randomly chosen code vector (or representative) of *p* dimensions. We then choose one observation randomly at each stage and compare it with all of the code vectors to determine the winning class, meaning the one whose code vector is closest (for a distance that has been determined beforehand). The code vectors of the winning class and of the neighbouring classes are moved in the direction of the chosen observation, so that the distance between them decreases.

This algorithm is analogous to a moving centres algorithm (in its stochastic version). However, the latter does not seek to conceptualise neighbourhoods of classes. Moreover, the only thing that it modifies at each stage is the code vector (or representative) of the winning class.

Following on from this, we assume that our readers are familiar with this algorithm (see inter alia Cottrell Fort, Pagès 1998).

Given that an arbitrary number of classes is chosen (it is often high because we frequently select grids of 8 by 8 or 10 by 10), we can reduce the number of classes, regrouping them by subjecting the vector codes to a classical hierarchical classification. We can then colour the class groups (called super classes) to enhance their visibility. Generally we observe that the only classes that such super classes regroup are contiguous ones. This can be explained by one of their properties, i.e., by the fact that the Kohonen algorithm respects the topology of the data space. Moreover, non-compliance with this property would indicate the algorithm's lack of convergence, or else a structure that has been particularly "folded" into the data set.

To describe the super classes, we calculate the basic statistics of the quantitative variables that are being used. We then study the way in which the modalities of the qualitative variables that the Kohonen classification algorithm does not use are distributed along the grid (Cottrell, Rousset 1997).

## 1.2 Classification of the observations that are being described by the qualitative variables - the KDISJ algorithm

This involves simultaneously classifying both individuals and the modalities of the qualitative variables that describe them. Analysts should be aware however that most of the time qualitative variables cannot be used in their existing form, even when the modalities are number coded. If no ordered relationship exists between the codes (for instance, 1 for blue eyes, 2 for brown eyes, etc.), it is no use applying them as if they were numerical variables, in a blind attempt to use Kohonen learning. Even if the codes were to correspond to an increasing or decreasing progression, this would only be meaningful if a linear scale were used (modality 2 corresponding to half of the progression between modalities 1 and 3). A fruitful method would then consist of processing the qualitative variables beforehand via a multiple correspondence analysis and preserving all of the co-ordinates. This is tantamount to coding all of the individuals by the co-ordinates that have been attributed to them as a result of this transformation. Once individuals have been represented by numerical variables, they can be classified using the Kohonen algorithm. We will however have lost the modalities, and the calculations will be both cumbersome and also costly in terms of calculating times - exactly that which we are trying to avoid by using the Kohonen algorithm.

The present paper introduces a method that has been adapted to qualitative variables, and which also enables a simultaneous processing of individuals and of modalities.

Consider $N$ individuals and a certain number $K$ of qualitative variables. Each variable $k = 1, 2, \ldots, K$ has $m_k$ modalities. Each individual chooses one and only one modality for each variable. If $M$ is the total number of modalities, each individual is represented by a $M$-vector comprised of 0s and 1s. There is only one 1 amongst the $m_1$ first components, only one 1 between the $(m_1+1)^{th}$ and the $(m_1+m_2)^{th}$, etc. The table with $N$ rows and $M$ columns that is formed in this way is the *complete disjunctive table*, called $D$. Note that it contains all of the information that will enable us to include individuals as well as the modalities' distribution.

We note $d_{ij}$ as the general term of this table. This can be equated to a contingency table that crosses an "individual" variable with $N$ modalities and a "modality" variable with $M$ modalities. The term $d_{ij}$ takes its values in $\{0,1\}$.

We use an adaptation of an algorithm (KORRESP) that has been introduced to analyse contingency tables which cross two qualitative variables. This algorithm is a very fast and efficient way of analysing the relationships between two qualitative variables. Please refer inter alia to Cottrell, Letrémy, Roy (1993) to see the various ways it can be applied to real data.

We calculate the row sums and column sums by:

$$d_{i.} = \sum_{j=1}^{M} d_{ij} \text{ et } d_{.j} = \sum_{i=1}^{N} d_{ij}.$$

Note that with a complete disjunctive table, $d_{i.}$ is equal to $K$, regardless of $i$. The term $d_{.j}$ represents the number of persons who are associated with the modality $j$.

In order to use a $\chi^2$-distance along the rows as well as down the columns, and to weight the modalities proportionately to the size of each sample, we adjust the complete disjunctive table, and put:

$$d_{ij}^c = \frac{d_{ij}}{\sqrt{d_{i.}d_{.j}}} \quad d_{ij}^c = \frac{d_{ij}}{\sqrt{d_{i.}d_{.j}}}.$$

When adjusted thusly, the table is called $D^c$ (adjusted disjunctive table). This transformation is the same as the one that Ibbou proposes in his thesis (Ibbou 1998; Cottrell, Ibbou, 1995).

These adjustments are exactly the same as the ones that correspondence analysis entails. This is in fact a principal weighted component analysis that uses the Chi-Square distance simultaneously along the row and column profiles. It is the equivalent of a principal components analysis of data that has been adjusted in this way.

We then choose a Kohonen network, and associate with each unit a code vector that is comprised of $(M + N)$ components, with the $M$ first components evolving in the space for individuals (represented by the rows of $D^c$) and the $N$ final components in the space for modalities (represented by the columns of $D^c$). The Kohonen algorithm lends itself to a double learning process. At each stage, we alternatively draw a $D^c$ row (i. e. , an individual), or a column (i. e. , a modality).

When we draw an individual $i$, we associate a modality $j(i)$, thus maximising the coefficient $d_{ij}^c$, i.e., the rarest modality out of all of the corresponding ones in the total population. We then create an extended individual vector of dimension $(M + N)$. Subsequently, we try to discover which is the closest of all the code vectors, in terms of the Euclidean distance (restricted to the $M$ first components). Note $u$ the winning unit. Next we move the code vector of the unit $u$ and its neighbours closer to the extended vector $(i, j(i))$, as per the customary Kohonen law.

When we draw a modality $j$ with dimension $N$, we do not associate an individual with it. Indeed, by construction, there are many equally placed individuals, and this would be an arbitrary choice. We then seek the code vector that is the closest, in terms of the Euclidean distance (restricted to the $N$ last components). We then move the $N$ last components of the winning code vector and its neighbours closer to the corresponding components of the modality vector $j$, without modifying the $M$ first components.

By so doing, we are carrying out a classical Kohonen classification of individuals, plus a classification of modalities, maintaining all the while their association with one another. After the convergence, the individuals and the modalities are classified into Kohonen classes. "Neighbouring" individuals or modalities are classified in the same class or in neighbouring classes. We call the algorithm that has been defined thusly KDISJ.

When we are not trying to classify individuals but only modalities, we can use another algorithm that draws its inspiration from the genuine Kohonen algorithm. This is called KMCA. We can then classify individuals as if they were additional data (for definitions and applications, see inter alia Ibbou's thesis, Ibbou, 1998). We can also classify individuals alone, and then classify as additional data the "virtual individuals" associated with the modalities that have been calculated from the rows of the Burt matrix. Finally we can classify modalities alone (as is the case with KMCA) and classify individuals subsequently, once they have been properly normalised. This is what Ibbou called KMCA1 and KMCA2. These methods generate findings that are very comparable to those that can be found with KDISJ, but they do require a few more iterations.

## 2. THE CLASSIFICATION

### 2.1 Classification using a Kohonen matrix and a regrouping into 10 super classes

On the map below (a 7 by 7 grid) we display findings from a simultaneous classification of individuals and variables. To simplify this representation, we have in each case displayed the current modalities, the number of individuals who have been classified, and between brackets the number of persons working on an OEC or FTC basis.

| AGE1 34(12,22) | 3(3,0) | DET4 34(31,3) | 5(5,0) | NITE2 32(29,3) | 2(2,0) | HORPOS 28(26,2) |
|---|---|---|---|---|---|---|
| 10(5,5) | 4(4,0) | 7(5,2) | 3(3,0) | 13(12,1) | SUN2 34(33,1) | 0 |
| NITE3 12(9,3) | 2(2,0) | ABS3 22(22,0) | 10(10,0) | DET2 38(26,12) | 0 | LEND2 36(33,0) |
| SUN3 19(18,1) | 1(1,0) | RECUP2 **31(31,0)** | 4(4,0) | 13(13,0) | 16(16,0) | 2(2,0) |
| 11(5,6) | SAT3 0 | FEM INVOL HORVAR JWK1 ABS1 9(8,1) | ABS2 2(24,1) | 8(8,0) | 44(43,1) | |
| 18(0,18) | 1(0,1) | 4(4,0) | OEC DET1 AGE3 LEND1 NITE1 RECUP1 SAT2 WED3 21(21,0) | AGE2 JWK2 WED2 29(29,0) | SAT1 INVOL 2(2,0) | DET3 31(31,0) |
| FTC 31(0,31) | 15(0,15) | MAN 38(38,0) | 1(1,0) | AGE4 SUN1 37(37,0) | HORIDE RECUP0 25(24,1) | WED1 30(26,4) |

*Figure 2.* Distribution of modalities and individuals across the grid

Note: The squares in gray feature a much higher percentage of OEC than the total population does.

Note how modalities and individuals are distributed amongst the various classes in a relatively balanced fashion. Fixed term contracts are mostly found to the left of the map. Remember that they only represent 17% of all contracts.

The modalities that correspond to the best working conditions (in other words, and for the purposes of the present paper, to more regular working times; to no night-time, Saturday or Sunday shifts; to open-ended contracts; and to voluntary part-time status) are associated with all age brackets, except for young persons, and are found to the bottom right. These correspond to relatively favourable work situations. Inversely, the young persons modality is located to the top right, and is associated with "unpleasant" modalities such as night shifts, Sunday shifts, no chance to take any time off etc.

The modality for women (who are present everywhere and who constitute the vast majority of the total population, to wit 88%) is close to the centre of the map and associated with the involuntary part-time modality that is close to the FTC modality.

### 2.2 Regrouping the classes

Next we diminish the number of classes by carrying out a hierarchical classification of the 49 code vectors. After several attempts to obtain a reasonably small number of classes, we have kept the 10 super classes that are represented below (see Figure 3).

*Figure 3.* The 10 super classes

The total population is relatively well balanced amongst these 10 super classes, with class 4 alone featuring a much larger sample. This will become understandable once we explain why – such individuals' working conditions are the most standard.

*Table 3.* Absolute frequencies

| Class | 1 | 2 | 3 | 4 | 5 | 6 | 7 | 8 | 9 | 10 |
|---|---|---|---|---|---|---|---|---|---|---|
| Size | 101 | 108 | 87 | 241 | 51 | 38 | 43 | 89 | 41 | 28 |

*Table 4.* Description of classes using qualitative variables (frequencies expressed as the percentage that the modality accounts for in each class)

|  | 1 | 2 | 3 | 4 | 5 | 6 | 7 | 8 | 9 | 10 | Tot |
|---|---|---|---|---|---|---|---|---|---|---|---|
| OEC | **99** | *40* | **100** | 92 | *47* | 92 | 77 | **94** | 88 | 93 | 83 |
| FTC | *1* | **60** | *0* | 8 | **53** | 8 | 23 | *6* | 12 | 7 | 17 |
| MAN | *0* | **54** | *2* | *1* | 14 | 13 | 16 | 10 | 12 | 7 | 12 |
| FEM | **100** | *46* | **98** | **99** | 86 | 87 | 84 | 90 | 88 | 93 | 88 |
| AGE1 | *0* | *0* | *0* | *0* | **100** | *0* | 5 | *0* | *0* | *0* | 6 |
| AGE2 | 36 | 51 | 43 | 45 | *0* | 26 | **65** | 39 | 39 | 39 | 40 |
| AGE3 | 42 | 22 | 36 | 32 | *0* | 39 | *9* | 44 | 29 | 43 | 31 |
| AGE4 | 23 | 27 | 22 | 22 | *0* | 34 | 21 | 17 | 32 | 18 | 22 |
| HORIDE | 52 | 61 | 48 | **75** | 35 | *26* | 49 | *29* | *29* | *0* | 52 |
| HORPOS | *0* | *0* | *0* | *0* | *0* | *0* | 12 | *0* | 8 | **100** | 4 |
| HORVAR | 48 | 39 | 52 | 25 | **65** | **74** | 40 | **71** | **71** | *0* | 44 |
| JWK1 | **96** | 83 | **89** | **91** | 78 | 79 | *47* | *39* | 68 | *57* | 79 |
| JWK2 | *4* | 17 | 11 | *9* | 22 | 21 | **53** | **61** | 32 | **43** | 21 |
| NITE1 | **100** | 92 | 95 | **99** | 86 | 95 | *51* | *64* | **100** | 75 | 90 |
| NITE2 | *0* | 5 | 4 | *1* | 10 | 5 | **21** | **27** | *0* | **21** | 7 |
| NITE3 | *0* | 3 | 1 | *0* | 4 | 0 | **28** | 9 | *0* | 4 | 3 |
| SAT1 | **88** | 57 | *3* | **80** | 29 | 53 | *5* | *2* | 34 | *21* | 49 |
| SAT2 | *6* | 23 | *5* | 10 | 14 | 26 | *5* | **91** | 27 | **68** | 23 |
| SAT3 | *6* | 20 | **92** | *10* | **57** | 21 | **90** | *7* | 39 | *11* | 28 |
| SUN1 | **96** | 88 | 82 | **99** | 69 | 87 | *0* | *13* | 83 | *57* | 76 |
| SUN2 | *4* | 12 | 18 | *1* | 18 | 10 | *2* | **87** | 17 | **43** | 18 |
| SUN3 | *0* | *0* | *0* | *0* | 13 | 3 | **98** | *0* | *0* | *0* | 6 |
| WED1 | **41** | 13 | 21 | **33** | 10 | 16 | 14 | 11 | 15 | 14 | 23 |
| WED2 | 14 | 9 | 10 | *9* | 18 | *8* | 12 | **46** | 19 | **43** | 16 |
| WED3 | 45 | **78** | 69 | 58 | 72 | **76** | 74 | 43 | 66 | 43 | 61 |
| ABS1 | 70 | 81 | 72 | 71 | 67 | 82 | *58* | 77 | 73 | 75 | 73 |
| ABS2 | **24** | 9 | *0* | **21** | 6 | 13 | 16 | 8 | 10 | 18 | 14 |
| ABS3 | *6* | 10 | **28** | 8 | **27** | *5* | **26** | 15 | 17 | 7 | 13 |
| DET1 | *0* | 81 | 90 | 75 | **88** | *37* | 72 | 68 | *0* | 78 | 63 |
| DET2 | *0* | 5 | *0* | **25** | 2 | 8 | 14 | 12 | *0* | 7 | 11 |
| DET3 | **100** | 14 | 10 | *0* | *4* | **45** | *5* | 13 | *0* | 11 | 19 |
| DET4 | *0* | *0* | *0* | *0* | 6 | 11 | 9 | 6 | **100** | 3 | 7 |
| INVOL | *12* | **80** | **74** | 44 | **82** | 53 | 60 | 34 | 46 | *21* | 50 |
| INVOL | **88** | *20* | *26* | 56 | *18* | 47 | 40 | 66 | 54 | **79** | 50 |
| LEND1 | 100 | 100 | 100 | 100 | 100 | *5* | 95 | 98 | 100 | 100 | 95 |
| LEND2 | 0 | 0 | 0 | 0 | 0 | **95** | 5 | 2 | 0 | 0 | 5 |
| RECUP0 | 44 | 57 | 40 | 67 | 65 | 34 | 44 | 32 | 37 | 61 | 52 |
| RECUP1 | 34 | 20 | 31 | 19 | 23 | 29 | 40 | 43 | 39 | 25 | 28 |
| RECUP2 | 22 | 23 | 29 | 14 | 12 | 37 | 16 | 25 | 24 | 14 | 20 |

Note: The numbers written in **bold font** correspond to particularly high values and those in *italics* to particularly low values.

We can verify that in most cases, the modalities find themselves either within or else close to one of the classes where they have a significant role to play. We can control this by calculating each modality's deviation for each of the 10 super-classes[1]. As can be expected, such deviations are positive 85 % of the time.

We then study the 5 quantitative variables' average values across the 10 classes:

*Table 5.* Description of the 10 classes by their quantitative variables (averages)

| Variable | 1 | 2 | 3 | 4 | 5 | 6 | 7 | 8 | 9 | 10 | Total |
|---|---|---|---|---|---|---|---|---|---|---|---|
| DMIN | 27.1 | 23.7 | 24.4 | 25.4 | 21.8 | 24.3 | 22.1 | 23.2 | 24.4 | 24.9 | 24.5 |
| DMAX | 29.1 | 27.7 | 27.4 | 27.2 | 24.1 | 29.5 | 32.2 | 32.3 | 28.9 | 32.0 | 28.5 |
| DTHEO | 27.0 | 24.4 | 24.5 | 25.7 | 22.3 | 24.8 | 25.4 | 25.8 | 25.1 | 27.1 | 25.4 |
| HSUP | 0.69 | 1.8 | 3.45 | 0.95 | 1.16 | 1.82 | 2 |  | 1.48 | 1.78 | 1.36 | 1.51 |
| HPROL | 1.65 | 1.97 | 1.54 | 0.78 | 0.75 | 2.08 | 2.53 | 2 |  | 2.71 | 0.86 | 1.5 |

Note that classes 6, 7, 8 and 10 display significant disparities between minimum and maximum workweek durations. Fisher statistics corresponding to these 5 variables show that they are all discriminatory in nature.

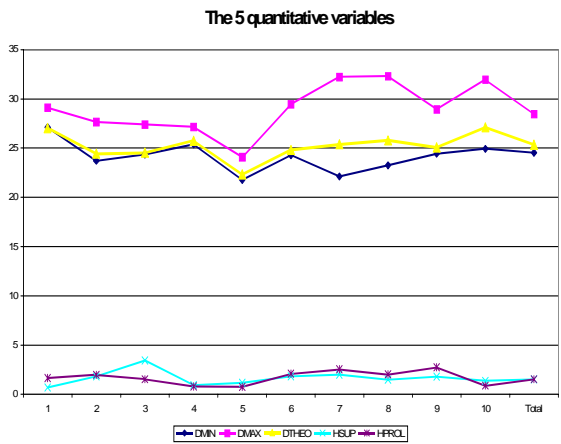

*Figure 4.* Quantitative variables in the 10 classes and in the total population

Based on these elementary statistics, it is possible to both describe the 10 classes and to develop a typology.

*Table 6.* Typology

| 1 | Employees, **choosing** voluntarily to work on a part-time basis; no Saturday shifts; they determine their own working schedules (very little overtime) |
|---|---|
| 2 | Men, working on a FTC basis; with Wednesday shifts; possibility of taking time off without any problem |
| 3 | Women who have had a part-time status **forced** upon them; every week they have the same number of workdays but daily schedules are variable; Saturday shifts; no possibility of taking any time off; no carryover of working hours (a lot of overtime) |
| 4 | The largest class, with 29% of the total. Employees working on an OEC basis; over the age of 25; no night-time or Sunday shifts; schedule is determined by the firm but flexibility is a possible; identical work schedules every week, but employees know their schedule for the next few days; time off can be taken under certain conditions; no reason to carry over working hours (shifts are rarely extended and there is little overtime) |
| 5 | All young persons under age of 25 (half OEC and half FTC); shift extensions are infrequent |
| 6 | Employees do not know their schedules for the next few days; average of almost 4 hours a week of shift extensions or overtime hours |
| 7 | It is **customary** for employees to work night and Sunday shifts (an average of more than 4h30 of shift extensions and overtime per week, some workweeks are more than 32h long even though they are allegedly doing part-time work). |
| 8 | Employees **sometimes** work night-time, Sunday, Saturday and Wednesday shifts, and do not work the same number of days every week (some weeks they can work more than 32h). |
| 9 | Working schedules are determined in a different way; shifts can be significantly extended |
| 10 | Everyone's schedule is **posted** openly (possibility of an approximately 32h workweek) |

On the super class representation, it is clear to see that the FTC and OEC modalities are distinct and separate (class 2 and class 4), as are men and women. As expected, women are associated with involuntary part-time work. The "voluntary part-time" modality can be found in class 1, near the OEC modality. Class 4 features the modalities that correspond to "normal" working conditions, with all ages being represented except for young persons.

For a more exhaustive summary, the reader can refer to the report edited by Cottrell, Letrémy, Macaire et al. (2001).

## CONCLUSION

In presenting our conclusions on part-time workers, we will refer to some of the descriptive statistics that the present paper was unable to mobilise, due to a lack of space.

Firstly, part-time work is more of an involuntary phenomenon for temporary employees than for permanent employees

The INSEE's Timetable survey raised a number of issues about part-time workers, during its attempt to test the "voluntary" nature of this form of employment. Regarding open-ended contracts, nearly 60% of all part-timers stated that this had been their choice, i. e., it was not imposed on them by their employer, either at the time of recruitment or else through the transformation of a full-time position into a part-time one. Around half stated that they had freely chosen their "shift system". In comparison, amongst employees working under fixed term contracts, fewer than 20% of all part-timers had volunteered for this status, but around 30% were working a shift system of their choice. This is quite a difference. Moreover, women unsurprisingly state more frequently than men do that they were the ones who had made the decision to work on a part-time basis. However, we know that choice is a highly relative concept, as all choice is made under constraint. We also know that family requirements often cause women to prefer part-time work.

Male and female part-time workers' situations vary greatly, depending on whether they are working under the aegis of an open-ended or a fixed term contract. Around 70% of women part-timers working on a fixed term contract do not get to choose their schedules (versus 54% of all men in this position). For part-timers working on an open-ended contract the gap is both reversed and smaller – in this population, 48% of women do not get to choose their schedules, versus 55% of men. It is as if the difference between OEC and FTC women were greater than between OEC and FTC men, whose situation is more homogeneous. Women on an open-ended contract basically choose their own schedule, more than men in this situation do. Working on a fixed term contract, however, they have less choice in their schedules.

In analysing the responses given to the question "Would you like to work more?", we learn that the more atypical the contract, the more employees would prefer to work more, as long as the increase in pay is proportional to the increase in the number of hours they work. Amongst atypical jobs, it is primarily part-time workers on fixed term contracts (and temporary workers, albeit to a lesser extent) who would like to work more. The same opposition between part-time and full-time work can be found in responses to questions relating to the desire to work less: unsurprisingly it is the part-timers who are less in favour of working fewer hours. Furthermore, those who are out looking for a new job are basically part-time employees on a fixed term contract (more than 40%) and temporary workers (more than 50%).

## APPENDIX

The data we used comes from the latest INSEE Timetable survey, the fourth of its kind (the previous one having been carried out in 1985-1986). It ran from February 1998 to February 1999 in 8 successive survey waves. Focusing on French lifestyle and working patterns, the full study looked at compensated professional working times, and more specifically at people's working times in their "*main current occupation*". The sample is comprised of the only salaried population that can provide comprehensive data on its professional working times. Teachers (who often make incoherent statements about their working times, equating them with contact hours alone) and other abnormal cases were taken out of the sample. 1,153 individuals were eliminated thusly, leaving a database of 5,558 wage-earning individuals.

When this sample is linked to data from the INSEE's 1998[2] or 1999 Employment surveys, no major difference is detected between the two in percentage terms. If we structure the data according to the type of work (full-time OEC, part-time OEC, full-time FTC, part-time FTC, temporary workers, other), we come up with two very similar distributions (see table 6 below).

*Table 7.* Distribution of sample according to form of employment in the 1998 INSEE Timetable and Job surveys

|  | *OEC FT* | *OEC PT* | *FTC FT* | *FTC PT* | *Temp* | *Others* |
|---|---|---|---|---|---|---|
| Timetable Survey | 4,033 | 690 | 258 | 137 | 115 | 325 |
| % of total sample | 72.6 | 12.4 | 4.6 | 2.5 | 2.1 | 5.8 |
|  | 85% | 7.1% | 2.1% | 5.8% |  |  |
| Job survey | 88.12%[3] | 5.57%[4] | 2.08% | 4.22% |  |  |

The breakdown between permanent/non permanent workers or between part-timers/full-timers is very comparable in the two surveys. Men represent a share of between 53 et 54% in both studies (and women between 45 and 46%). However, in terms of respondents' ages, the Timetable survey slightly over-represents people between the age of 40 and 49 (by 3 points) and under-represents the 25-39 age bracket.

Other differences can be observed:
− over-representation of the industrial sector (by 6 points) in the Timetable survey;
− under-representation by 5 points of the service sector;
− under-representation by around 5 points of sectors such as healthcare, education and social work.

## NOTES

1. The deviation for a modality *m* (shared by $n_m$ individuals) and for a class *k* (with $n^k$, individuals) can be calculated as the difference between the number of individuals who possess this modality and belong to the class k and the "theoretical" number $n_m n_k / n$ which would correspond to a distribution of the modality *m* in the class *k* that matches its distribution throughout the total population.
2. Source: Employment Survey 1998, INSEE findings, n° 141-142, 1998, 197 pages.
3. Except for non-tenured State and local authority employees.
4. TC except for State and local authority officials, + non-tenured State and local authority employees.